\documentclass[11pt]{article}
\usepackage{graphicx}
\usepackage{color}
\usepackage{amsmath,amsfonts,amssymb}
\def\la{\big\langle}
\def\ra{\big\rangle}

\def\ds{\displaystyle}
\def\forall{\hbox{for all}~}
\def\L{{\bf L}}

\def\bfv{{\bf v}}

\def\bfw{{\bf w}}
\def\bfu{{\bf u}}
\def\bfb{{\bf b}}

\def\bfn{{\bf n}}
\def\bfe{{\bf e}}

\def\bfp{{\bf p}}

\def\bfU{{\bf U}}

\def\bpm{\begin{pmatrix}}
\def\epm{\end{pmatrix}}

\def\S{{\cal S}}
\def\F{{\cal F}}

\def\M{{\cal M}}
\def\N{{\cal N}}

\def\R{{\mathbb R}}

\def\rank{\hbox{\rm Rank}}
\def\implies{\Longrightarrow}

\def\vs{\vskip 2em}

\def\v{\vskip 1em}

\def\begi{\begin{itemize}}
\def\endi{\end{itemize}}
\def\C{{\cal C}}
\def\div{\hbox{div}}
\def\ov{\overline}

\def\bega{\begin{array}}
\def\enda{\end{array}}

\def\bel{\begin{equation}\label}
\def\eeq{\end{equation}}
\def\sqr#1#2{\vbox{\hrule height .#2pt
\hbox{\vrule width .#2pt height #1pt \kern #1pt
\vrule width .#2pt}\hrule height .#2pt }}
\def\square{\sqr74}
\def\endproof{\hphantom{MM}\hfill\llap{$\square$}\goodbreak}

\newtheorem{theorem}{Theorem}[section]

\newtheorem{lemma}{Lemma}[section]

\newtheorem{remark}{Remark}[section]

\begin{document}
\title{\bf  Generic Singularities for 2D Pressureless Flow}
\vs

\author{Alberto Bressan$^{*}$, Geng Chen$^{**}$, and Shoujun Huang$^{\dag}$\\
\, \\
$^*$Department of Mathematics, Penn State University, \\
University Park, PA ~16802, USA.\\
$^{**}$Department of Mathematics, University of Kansas,\\
 Lawrence, KS 66045, USA.\\
$^\dag$College of Mathematical Medicine, Zhejiang Normal University, \\
Jinhua 321004, P. R. China.\\
\, \\
E-mails: axb62@psu.edu,~gengchen@ku.edu,~sjhuang@zjnu.edu.cn.}
\maketitle
\v

{\bf Abstract:}  In this paper, we consider the Cauchy problem for pressureless gases in two space dimensions with generic smooth initial data (density and velocity).
These equations give rise to singular curves, where the mass has positive density w.r.t.~1-dimensional Hausdorff measure.
We observe that the system of equations describing these singular curves is not hyperbolic.
For analytic data, local solutions are constructed using a version of the Cauchy-Kovalevskaya theorem.
We then study the interaction of two singular curves, in generic position.
Finally, for a generic initial velocity field,  we investigate the asymptotic structure of the smooth solution up to the first time when a singularity is formed.

{\bf Key words:}  Pressureless gases; formation of singularities; generic data;
Cauchy-Kovalevskaya theorem.

\section{Introduction}
\label{s:1}
\setcounter{equation}{0}

We consider the initial value problem for the equations of pressureless gases
in two space dimensions:
\begin{equation}\label{SP}
\left\{\begin{array}{rl}
\!\!\partial_t\rho+\div( \rho \bfv)&=0, \cr
\!\!\partial_t(\rho \bfv)+\div\big( \rho \bfv\otimes \bfv\big)&=0,\enda
\right.  \qquad\quad  t\in  ]0,T[\,,\quad x\in\R^2,\eeq
\bel{ic}\rho(0,x)~=~\bar \rho(x), \qquad \bfv(0,x)~=~\bfw(x).\eeq

 A simple  measure-valued solution of (\ref{SP})
is provided by a finite collection of particles, moving with constant speed
in the
absence of forces. Whenever two or more particles collide,
they stick to each other as a single compound particle. The mass of the new particle
is equal to the sum of the masses of the particles involved in the collision, while
its velocity is determined by the conservation of momentum.
For this reason, one often refers to
(\ref{SP}) as the ``sticky particle" system.

A general Lagrangian framework for solving these equations was introduced in
\cite{Sever}.  However, as later shown in
\cite{BN}, for the Cauchy problem (\ref{SP})-(\ref{ic}) both the existence and the uniqueness of solutions can fail, even for $\L^\infty$ initial data.   A simplified version
of the counterexamples introduced in \cite{BN}  is shown in Fig.~\ref{f:z112}.

\begin{figure}[ht]
\centerline{\hbox{\includegraphics[width=16cm]
{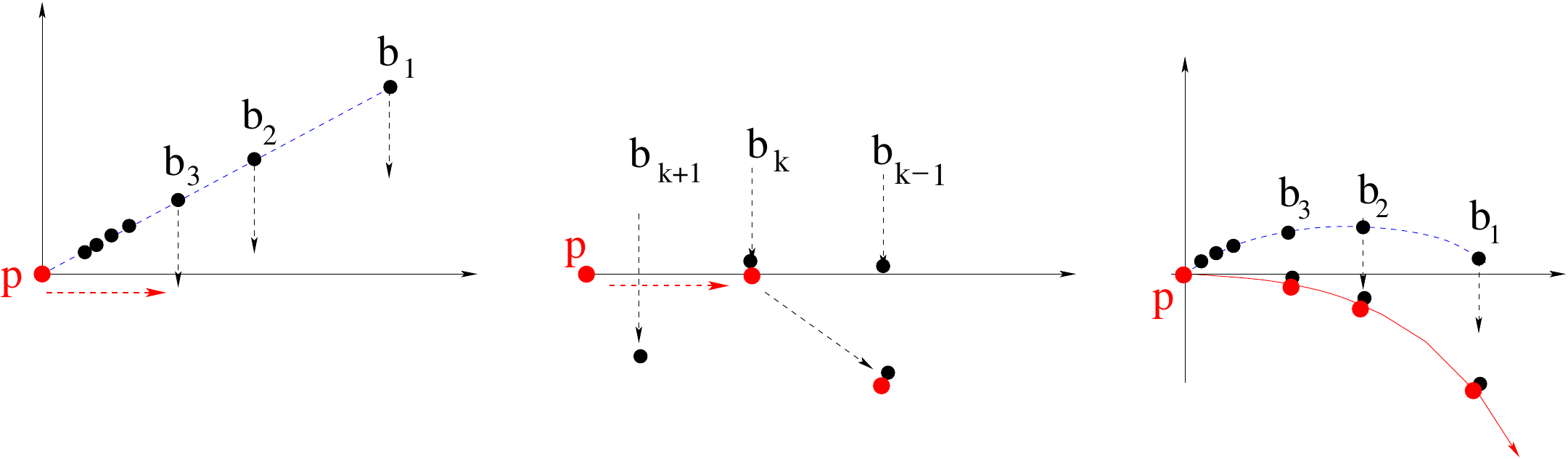}}}
\caption{\small Left: an example of initial data without any solution. At the initial time $t=0$
a ``target" particle $p$ with unit mass
is located at the origin, with speed $\bfv_0=(1,0)$.   A countable set of ``bullet" particles $b_k$, with masses
$m_k= k^{-2}$ are initially located at points $(k^{-1}, k^{-1})$, with speed $\bfv_k=(0, -1)$.
As shown in the center figure, the bullet $b_k$ hits the target  at time $t_k=k^{-1}$
if and only if none of the bullets $b_j$ with $j>k$ has already scored a hit.
Since there is no ``first bullet" that can hit the target particle, every hypothetical
interaction pattern leads to a contradiction. Right: an example of initial data
with two solutions.   Here the initial position of the bullet particles $b_k$ is chosen so that
they can all  hit the target, one after the other.   Eventually, all particles thus stick to each other
as a single
point mass, traveling with the speed of the barycenter.
However, there exists also a second solution where
all particles proceed in rectilinear motion, without ever touching each other.}
\label{f:z112}
\end{figure}

When the initial data range in a space of measures, the recent
paper \cite{BD} shows that, in dimension $n\geq 2$, generic solutions
correspond to free flow.  In other words, for ``almost all" measure-valued initial data,
particles never hit each other and
their motion is rectilinear, with constant speed.

In the present note we also consider generic initial data, but at the opposite end of the regularity spectrum.
Namely, we shall assume that  $\bar \rho, \bfw$ are smooth, or even analytic functions.
As long as the solution remains smooth, it can be directly computed by the method
of characteristics.  Our main interest is the local construction of solutions
beyond singularity formation.

In one space dimension,
the solution will generically contain a finite number of point masses,
which can eventually collide with each other \cite{D}.
In two space dimensions, the equations (\ref{SP}) give rise to
``singular curves", where the mass has positive density w.r.t.~1-dimensional
Hausdorff measure.

Our goal is to understand the formation
and the time evolution of these singular curves.
A key feature, observed in the 2-dimensional case, is that
the equations describing these singular curves are not of hyperbolic type.
Therefore, apart from some special configurations, the Cauchy problem is ill posed.
Local solutions will thus  be constructed only within the class of analytic data.
In a way, this situation resembles the case
of vortex sheets \cite{CO}.

As usual, the space of three times continuously differentiable vector fields
$\bfw=(w_1,w_2)$ on $\R^2$  will be denoted by
$$\C^3(\R^2;\R^2)~\doteq~\left\{ \bfw:\R^2\mapsto\R^2\,;~~\|\bfw\|_{\C^3}\doteq
\sup_{x\in \R^2} \sup_{|\alpha|\leq 3, ~i=1,2} |D^\alpha w_i(x)|~<~+\infty\right\}.$$
It will be convenient to work within the subspace of $\C^3$
vector fields whose Jacobian
matrix $D\bfw$ vanishes as $|x|\to +\infty$.  This space will be denoted by
\bel{C30}\C^3_0(\R^2;\R^2)~\doteq~\left\{ \bfw\in \C^3(\R^2;\R^2),~~\lim_{|x|\to +\infty} D\bfw(x)~=~ {\bf 0}
\right\}.\eeq
Note that $\C^3_0$ is a closed subspace of $\C^3$.  Hence it is a Banach space as well.
The advantage of working in this smaller space is that, for any given
$\bfw\in \C^3_0$ and $T>0$, we can find $R>0$ large enough so that
for all $t\in [0,T]$ one has
\bel{wreg}\det\bigl(I+tD\bfw(x)\bigr)~>~0\qquad\hbox{for all}~~t\in [0,T],~~|x|\geq R.\eeq
As a consequence, for $t\in [0,T]$, the only singularities of the solution to
(\ref{SP})-(\ref{ic}) can occur within a bounded domain.
\v
The main issues we seek  to  understand are:
\begi
\item[(i)] Time evolution of a singular curve.
\item[(ii)] Interaction of two singular curves.
\item[(iii)] Asymptotic structure of the solution in a neighborhood of
a point $(\tau,y)\in \R_+\times\R^2$
where the $2\times 2$ Jacobian matrix $D_x\bfv$ blows up and a new singularity is formed.
\item[(iv)] Structure of a singular curve immediately after its formation.
\endi
In the present paper we focus on (i)--(iii). The initial stages of new singular curve,
requiring a deeper analysis, will
be described  in the forthcoming paper \cite{BCH}.
In Section~\ref{s:2} we briefly review the construction of smooth solutions,
by the method of characteristics.
In Section~\ref{s:3} we derive the system of equations describing a singular curve.
We observe that in general this system is not hyperbolic. Because of this, local solutions will
be constructed for analytic initial data, relying on a version of the Cauchy-Kovalevskaya
theorem~\cite{evans, RR}.  Section~\ref{s:4} deals with
the interaction of two singular curves, in two main
cases: at the first time where the curves  touch each other tangentially,
and at a later time when part of the two curves have already merged into a single one.
These are illustrated in Figures~\ref{f:gen38} and \ref{f:gen39}, respectively.

Finally, in
Section~\ref{s:5}, we study the asymptotic structure of a smooth solution up to the first
time when a singularity is formed.
This happens along a characteristic
starting at a point $\bar x$
where the smaller eigenvalue of the Jacobian matrix
$D\bfw(\bar x)$ has a negative minimum $\lambda_1(\bar x)<0$.
Relying on techniques from differential geometry \cite{Bloom, GG},
we describe the behavior of a generic vector field $\bfw$
in a neighborhood of such a point $\bar x$.
Here ``generic"  means that the results apply to an open dense set of
vector fields in $ \C^3_0(\R^2;\,\R^2)$.
 We remark that, while the construction of a singular curve  requires
 analytic data,
the main results in Section~\ref{s:5} remain valid
more generally for initial data $\bfw\in \C^3$.

The present results can be seen as part of a general research program aimed at understanding
the singularities of solutions to nonlinear PDEs, for generic smooth initial data~\cite{Damon, Guck}.
In one space dimension,
singular solutions to the sticky particle model were studied in \cite{D}, while the recent paper
\cite{AR} derives the evolution equations for a singular surface in multidimensional space.
The papers \cite{BC, BHY} deal with generic singularities of solutions to a second order
variational wave equation.
Generic solutions to 1-dimensional  conservation laws were first considered
in \cite{GS, Gk2, S}.
In one space dimension, the analysis in \cite{Kong} provided  an asymptotic description
of this blow up, and also of the initial stages of the new shock, in a generic setting.
A study of shock formation for multidimensional scalar conservation laws can be found in
\cite{DM}.
For the 3-dimensional equations of isentropic gas dynamics,
the asymptotic of blow up, leading to shock formation, was studied in \cite{BSV}.
For earlier results on blowup for hyperbolic equations we refer to~\cite{A}.

The equations of pressureless gases have found applications as models
for the motion of granular media \cite{ERS}.
In particular, they provide
an approximate description of the large-scale distribution of matter in the universe
\cite{GSS, Z}.
More recently,  they have been shown to describe a limit of collective dynamics with short-range interactions \cite{ST}.

\section{Smooth solutions}
\label{s:2}
\setcounter{equation}{0}

We review the construction of a smooth solution to (\ref{SP}), by the method of characteristics.
Call $\bfv=(t,x)$ the velocity function.
Let $t\mapsto x(t)$ be a characteristic curve, such that
\bel{ch1}
{d\over dt} x(t)~=~\bfv(t, x(t)).\eeq
Assuming that the solution remains smooth, from (\ref{SP}) it follows
\bel{ch2}
{d\over dt} \rho(t, x(t))~=~\rho_t + \bfv\cdot \nabla\rho~=\,-\div(\rho \bfv)+\bfv\cdot \nabla\rho~=\,-\rho\, \div\, \bfv.\eeq
\bel{ch3}
{d\over dt}(\rho\bfv)(t, x(t))~=~(\rho\bfv)_t +\bfv\cdot D(\rho\bfv)~=~
-\div(\rho\bfv\otimes \bfv) +\bfv\cdot D(\rho\bfv)~=~-\rho\bfv\,\div\,\bfv.\eeq
\bel{ch4}
{d\over dt} \bfv(t, x(t))~=~{d\over dt}{(\rho \,\bfv)(t, x(t))\over \rho(t, x(t))}~=~
{(-\rho\bfv\, \div\,\bfv)\rho - (\rho \bfv)(-\rho\,\div\, \bfv)\over \rho^2}~=~0.
\eeq
By (\ref{ch4}) it follows
\bel{ch6} \bfv(x+ t\bfw(x))~=~\bfw(x).\eeq
In addition, conservation of mass implies
\bel{ch7}\rho(x+t\bfw(x))~=~{\bar \rho(x)\over \det\bigl(I + tD\bfw(x)\bigr)}\,.\eeq
The equations (\ref{ch6})-(\ref{ch7}) describe the solution up to the first time
$\tau$
where a singularity occurs.    This happens as soon as the matrix
$I+ tD\bfw(x)$ is no longer invertible.
This blow up time $\tau$ can thus be characterized as the smallest time $t>0$
such that $-1/t$ is an eigenvalue of $D\bfw(x)$, for some $x\in\R^n$.

\section{Evolution of a singular curve}
\label{s:3}
\setcounter{equation}{0}
We start by deriving  a system of equations for the evolution of a singular curve
$\gamma$ in 2-dimensional space, as shown in
Fig.~\ref{f:gen37}.  An alternative derivation has recently appeared also in \cite{AR}.
At time $t$
 the curve will be parameterized as
$$\xi~\mapsto ~y(t,\xi),\qquad\qquad \xi\in [a(t), b(t)],$$
so that, for each fixed $\xi$, the map  $t\mapsto y(t,\xi)$ traces a particle
trajectory.
We denote by $\eta=\eta(t,\xi)$
the linear density and by $\bfv(t,\xi)$ the velocity of particles along the curve.
This implies
\bel{91}
y_t(t,\xi)~=~\bfv(t,\xi).\eeq
Moreover, at time $t$, the total mass contained  on a portion of the
curve is
\bel{mm}m\Big(\{y(t,\xi)\,;~\xi\in [\xi_1, \xi_2]\}\Big)~=~\int_{\xi_1}^{\xi_2} \eta(t,\xi)\, d\xi.\eeq

To derive an evolution equation for  $\eta,\bfv$, let
$$\bfv^+(t,x),\quad \bfv^-(t,x),\qquad \rho^+(t,x),\quad \rho^-(t,x),$$
be respectively the velocity and the density (w.r.t.~2-dimensional Lebesgue measure)
of the gas, in front and behind the curve.

\begin{figure}[ht]
\centerline{\hbox{\includegraphics[width=7cm]
{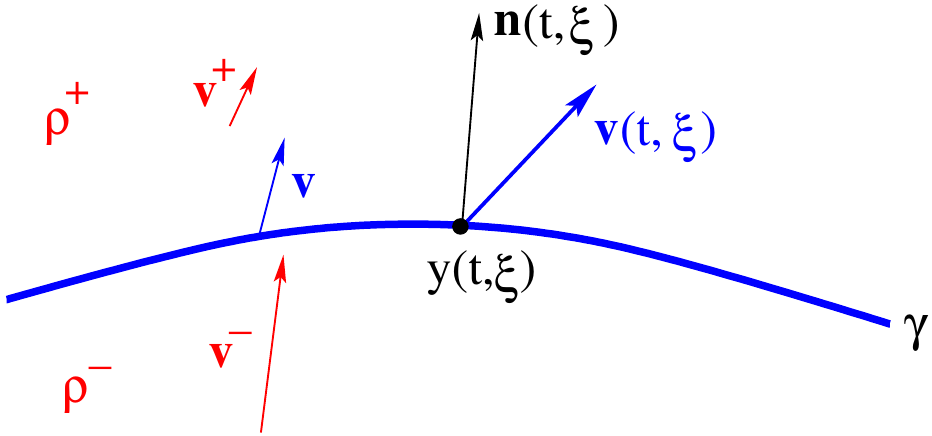}}}
\caption{\small  The variables in the system (\ref{ssy}), determining the
evolution of a singular curve $\gamma$.}
\label{f:gen37}
\end{figure}

Calling $\bfn(t,\xi)$ the unit normal vector to the curve at the point $y(t,\xi)$,
as in Fig.~\ref{f:gen37},
we assume that the inner products satisfy the admissibility condition
\bel{inp}\Big\langle \bfn(t,\xi),\, \bfv^+(t, y(t,\xi))\Big\rangle ~\leq~\Big\langle  \bfn(t,\xi),\, \bfv(t,\xi)\Big\rangle
~\leq~\Big\langle  \bfn(t,\xi),\, \bfv^-(t, y(t,\xi))\Big\rangle .\eeq
These inequalities have a clear  physical meaning. Indeed,
they imply that particles impinge on the singular curve from both sides
(and then stick to it).

It will be convenient to write the next equation using the
cross product of two vectors $\bfu= (u_1, u_2)$, $\bfv = (v_1, v_2)$ in $\R^2$:
$$\bfu\times \bfv~=~u_1 v_2-u_2 v_1~=~\la \bfu^\perp,\bfv\ra,\qquad \bfu^\perp = (-u_2, u_1).$$
Moreover, we shall denote by
$d\sigma$ the arc-length along the singular curve $y(t,\cdot)$, so that
$$ \bfn(t,\xi)\,d\sigma~=~y_\xi^\perp(t,\xi)\,d\xi\,.$$
Fix $\xi_1<\xi_2$ and consider the
portion of the singular curve  $\{y(t,\xi)\,;~\xi\in[\xi_1,\xi_2]\}$.
 Recalling (\ref{mm}), let
 $$M(t)~=~\int_{\xi_1}^{\xi_2} \eta(t,\xi)\, d\xi$$
be the total mass contained in this portion of the curve at time $t$.

Under the assumption (\ref{inp}), by conservation of mass, the rate at which
$M(t)$  increases is computed by
$$\bega{rl}\ds {d\over dt} M(t)&=~\ds\int_{\xi_1}^{\xi_2}\eta_t\, d\xi\\[4mm]
&=~\ds\int_{\xi_1}^{\xi_2}\la \bfv-\bfv^+, \bfn\ra\rho^+d\sigma
+\int_{\xi_1}^{\xi_2}\la \bfv^--\bfv, \bfn\ra\rho^-d\sigma
\\[4mm]
&=~\ds\int_{\xi_1}^{\xi_2}\la \bfv-\bfv^+, \, y_\xi^\perp\ra\rho^+d\xi+
\int_{\xi_1}^{\xi_2}\la \bfv^--\bfv, \, y_\xi^\perp\ra\rho^-d\xi
 \\[4mm]
&=~\ds
\int_{\xi_1}^{\xi_2}y_{\xi}\times (\bfv-\bfv^+)\rho^+d\xi+
\int_{\xi_1}^{\xi_2}y_{\xi}\times (\bfv^--\bfv)\rho^-d\xi \,.\enda$$
Since the above identity holds for all $\xi_1<\xi_2$, this implies
\bel{92}\eta_t(t,\xi)~=~[ y_\xi\times  (\bfv-\bfv^+)] \,\rho^+-
[ y_\xi\times  (\bfv-\bfv^-)] \,\rho^-,
\eeq
where the functions $\rho^\pm, \bfv^\pm$
 are evaluated at $(t,y(t,\xi))$.

A similar argument shows that, by conservation of momentum, one has
\bel{93}
\bigl[\eta(t,\xi)\, \bfv(t,\xi)\bigr]_t ~=~[ y_\xi\times  (\bfv-\bfv^+)] \,\rho^+\bfv^+-
[ y_\xi\times  (\bfv-\bfv^-)] \,\rho^-\bfv^-.\eeq
Combining (\ref{92})  with (\ref{93}), one obtains an
expression for  $\bfv_t(t,\xi)$, describing the acceleration of a particle along the singular curve.
\bel{94}\bega{rl} \ds \bfv_t(t,\xi)&\ds=~\left( {\eta\bfv\over\eta}\right)_t~=~{(\eta\bfv)_t\over\eta} - {\eta_t\over \eta}\, \bfv
\\[4mm] &\ds=~{1\over\eta} \bigg(
\bigl[ y_\xi\times  (\bfv-\bfv^+)\bigr]\, \rho^+(\bfv^+-\bfv)+
 \bigl[y_\xi\times  (\bfv-\bfv^-)\bigr]\, \rho^-(\bfv^--\bfv)\bigg).
\enda \eeq
Setting $\bfp=y_\xi$, and observing that $y_{\xi t}~=~\bfv_\xi$,
we thus obtain the semilinear system
\bel{ssy}
\left\{\bega{rl}y_t&=~\bfv,\\[3mm]
\eta_t&=~\bfp\times (\bfv-\bfv^+)\rho^+-\bfp\times (\bfv-\bfv^-)\rho^-,\\[3mm]
\bfv_t&=~\ds \Big(\bfp\times (\bfv-\bfv^+)\Big)\rho^+\, {\bfv^+-\bfv\over\eta}+\Big(\bfp\times (\bfv-\bfv^-)\Big)\rho^-\, {\bfv^--\bfv\over\eta}\,,\\[3mm]
\bfp_t - \bfv_\xi&= ~0.\enda\right.\eeq

We observe that the semilinear
system (\ref{ssy}) is not hyperbolic.
As a consequence,
in general we expect that  the Cauchy problem describing the evolution of a
singular curve $\gamma$  will be ill posed.
Existence and uniqueness of local solutions can be still achieved, but
 for analytic initial data.
For non-analytic initial data, the curve $\gamma$  might
instantly disintegrate.  In this case, in view of the examples in \cite{BN},
it is not even clear whether a weak solution exists at all.




Given  analytic initial data
\bel{icon}y(t_0,\xi)\,=\,y_0(\xi),\qquad \eta(t_0,\xi)\,=\,\eta_0(\xi),\qquad
\bfv(t_0,\xi)\,=\,\bfv_0(\xi),\eeq
a local solution  to the Cauchy problem (\ref{ssy})-(\ref{icon})
can be obtained by an application of the classical Cauchy-Kovalevskaya theorem.

Without loss of generality, we shall assume
$$t_0\,=\,0,\qquad y(t_0,0) \,=\, y_0(0)\,=\,0.$$
In addition, we assume that functions $\rho^\pm$
and $\bfv^\pm$ are given, so that the following conditions hold.
\begi
\item[{\bf (A1)}] {\it  The initial data  $y_0(\cdot)$,  $\eta_0(\cdot)$, $\bfv_0(\cdot)$ are analytic
in a neighborhood of $\xi=0$.  Moreover, $\eta_0(0)>0$. }
\item[{\bf (A2)}]  {\it  The functions
$\rho^{\pm}:\R\times\R^2\mapsto\R_+$ and $\bfv^{\pm}:\R\times\R^2\mapsto\R^2$
are analytic  in a neighborhood  of the point $(t,x)=(0,0)\in \R\times\R^2$.}
\end{itemize}
The main result of this section is
\begin{theorem}\label{t:31} Under the assumptions {\bf (A1)--(A2)}, the Cauchy problem (\ref{ssy})-(\ref{icon})
admits a real analytic solution defined in a neighborhood of the origin.
 \end{theorem}

{\bf Proof.} It suffices to rewrite the Cauchy problem in the standard form  to which the Cauchy-Kovalevskaya theorem applies.
Toward this goal,
it is convenient to call
$$x_1~=~\xi,\qquad x_2~=~t,$$ the independent variables,
and introduce the  new  dependent variables
$$\left\{\bega{rl}Y(t,\xi)&=~y(t,\xi)-y_0(\xi),\\[2mm] W(t,\xi)&=~\eta(t,\xi)-\eta_0(\xi),\\[2mm]
Z(t,\xi)&=~\bfv(t,\xi)-\bfv_0(\xi),\enda\right. \qquad \qquad \left\{\bega{rl}
{\bf P}(t,\xi)&=~\bfp(t,\xi)- y_{0,\xi}(\xi)\,,\\[3mm] \tau(t,\xi) &= ~t\,.\enda\right.$$
The above definitions ensure that  the initial data at $x_2=t=0$ vanish:
\bel{ida2} Y(x_1,0)\,=\,Z(x_1,0)\,=\,{\bf P}(x_1,0)\,=\,0\,\in\,\R^2,\
\quad \qquad W(x_1,0)\,=\,\tau(x_1,0)\,=\,0.\eeq
We now observe that, under the assumptions {\bf (A1)-(A2)}, the right hand side
of (\ref{ssy}) is an analytic function of the variables $t,y,\eta, \bfv, \bfp$,
because $\eta$ does not vanish near the origin.
Adding the trivial equation $\tau_t=1$, in the new variables the semilinear system
(\ref{ssy}) can be written more compactly as
\bel{cp}\left\{\begin{array}{rll}\bfU_{x_2}\!&=~{\bf B}\,\bfU_{x_1}+{\bf F}(x_1,\bfU)
\qquad &\hbox{for}~~x_2>0,\\[2mm]\bfU&=~{\bf 0}\qquad
&\hbox{for}~~x_2=0.\end{array}\right.\eeq
Here
$$\bfU~=~(Y,W,Z,{\bf P}, \tau)~\in~\R^2\times \R\times \R^2\times \R^2\times\R.$$
Moreover,
 ${\bf F}:\R\times\R^8\mapsto\R^8$
is an analytic function of all its arguments, while
 ${\bf B}= [b_{ij}]$ is a constant $8\times 8$ matrix. Namely, recalling the vector
 equation $\bfp_t-\bfv_\xi=0$ in (\ref{ssy}), one obtains
 $$ b_{6 4} \,=\, b_{75}\,=\,1,\qquad\qquad b_{ij}\,=\,0\qquad\hbox{otherwise}.$$
An application of the Cauchy-Kovalevskaya theorem \cite{evans, RR, W}
now yields the result.\endproof

\begin{remark} {\rm In connection with the original system (\ref{SP}), assume that the initial data $\bar\rho, \bfw$ in (\ref{ic}) is analytic.
Then the solution $(\rho,\bfv)$ constructed in Section~\ref{s:2} by the method of characteristics will remain analytic as long as the gradient remains bounded.
In particular, if $\rho^\pm, \bfv^\pm$ are solutions to (\ref{ch6})-(\ref{ch7}), then
the regularity assumption {\bf (A2)} will hold.

It is worth noting that analytic solutions to (\ref{ssy})  can be constructed even
if the admissibility conditions (\ref{inp}) fail.  In this case, however,
the solution would lose its physical meaning.  For this reason, in addition to the assumptions {\bf (A1)-(A2)}, it is natural to impose
\begi
\item[{\bf (A3)}]  {\it At the point $(t,x)=(0,0)$, one has }
\bel{inpa}\langle \bfn,\, \bfv^+\rangle ~<~\langle  \bfn,\, \bfv\rangle
~<~\langle  \bfn,\, \bfv^-\rangle .\eeq
\endi
}
\end{remark}

\section{Interaction of singular curves}
 \label{s:4}
\setcounter{equation}{0}
Having constructed a local solution which includes a singular curve, we now
study the interaction of two singular curves.

We consider two generic cases:
\begi
\item[(i)] Two curves touching each other at an interior point, tangentially
(Fig.~\ref{f:gen38}).
\item[(ii)] Two curves touching at an endpoint, with a nonzero angle (Fig.~\ref{f:gen39}).
\endi

The procedure for constructing local analytic solutions is the same in both cases.
We assume that, in a neighborhood of an interaction point, one is given
the analytic vector fields
 $\bfv^-, \bfv_0, \bfv^+$
and the analytic, positive density functions $\rho^-, \rho_0, \rho^+$.

This allows to construct separately the two singular curves
$\gamma_1=\gamma_1(t,\xi)$ and $\gamma_2=\gamma_2(t,s)$.
We define the functions $\xi\mapsto t^\sharp(\xi)$,
$\xi\mapsto s^\sharp(\xi)$, and the intersection curve $\gamma^\sharp$
implicitly by the equations
\bel{ssh}
\gamma_1(t^\sharp(\xi), \xi)~=~\gamma_2(t^\sharp(\xi), s^\sharp(\xi)),\eeq
\bel{gsh}
\gamma^\sharp(t^\sharp(\xi), \xi)~\doteq~\gamma_1(t^\sharp(\xi), \xi).\eeq
We think of  $t^\sharp(\xi)$ as the time when the particle $\xi$ on the first
curve hits some particle on the second curve.
In the following, we study two basic configurations where the equations
(\ref{ssh})-(\ref{gsh})
can be uniquely solved, thus determining the functions $t^\sharp, s^\sharp$.

\begin{figure}[ht]
\centerline{\hbox{\includegraphics[width=8cm]
{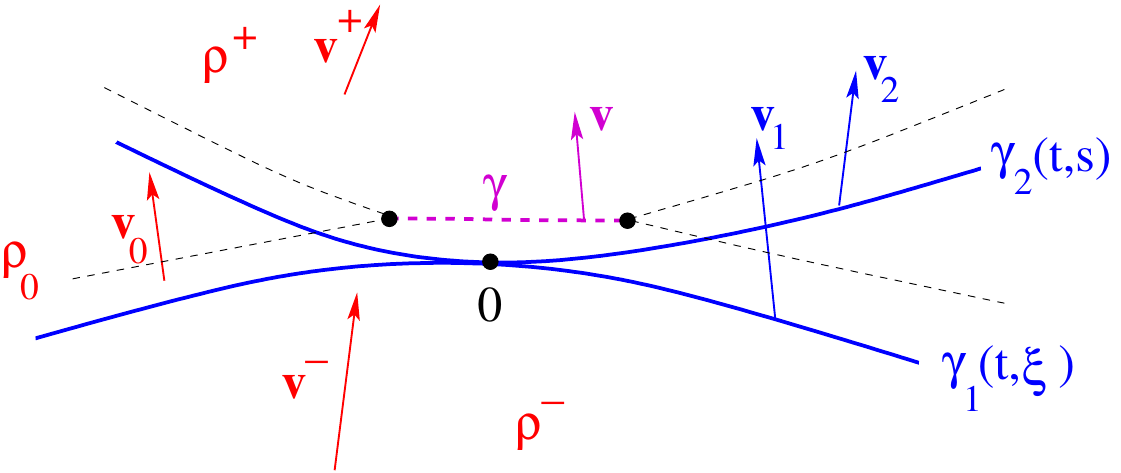}}}
\caption{\small  Two singular curves $\gamma_1,\gamma_2$,
which initially are tangent to each other at the origin. The dotted lines show the configuration at a later time $t>0$. Assuming that
the normal velocities satisfy (\ref{ndeg}), a new singular curve $\gamma$ is generated,
collecting particles from both $\gamma_1$ and $\gamma_2$, sticking to each other. }
\label{f:gen38}
\end{figure}

\subsection{Two singular curves tangent at one point.}
By a change of coordinates we assume that, at the initial time $t=t_0$, the curves $\gamma_1,\gamma_2$
can be parameterized respectively as
\bel{ttc}
\xi\mapsto \gamma_1(t_0,\xi)~=~\bigl(\xi, \phi_1(\xi)\bigr),\qquad\quad
s\mapsto \gamma_2(t_0,s)~=~\bigl(s, \phi_2(s)\bigr),\eeq
with
\bel{phi12}\phi_1(0)\,=\, \phi_2(0)\,=\,0, \qquad \phi_1'(0)\,=\, \phi_2'(0)\,=\,0,
\qquad \phi_1''(0)~<~\phi_2''(0).\eeq
Call $\bfn=(0,1)$ the unit vector perpendicular to both curves at the origin.
Assume that, always at the origin,  the velocities of the two curves at the initial
time $t=t_0$ satisfy
\bel{ndeg}\la \bfn,\,  \bfv_1\ra~>~\la \bfn,\,  \bfv_2\ra.\eeq
We claim that the vector equation (\ref{ssh}) uniquely determines
the two scalar functions $t^\sharp(\xi), s^\sharp(\xi)$, by the implicit function theorem, in a neighborhood of the origin.

Indeed, calling $\{\bfe_1,\bfe_2\}$ the standard orthonormal basis of $\R^2$, the
vector equation (\ref{ssh}) is equivalent to
\bel{Phi12}\left\{\bega{rl} \Phi_1&\doteq~\la\bfe_1, \gamma_1(t^\sharp, \xi)-\gamma_2(t^\sharp, s^\sharp)\ra~=~0,\\[2mm]
 \Phi_2&\doteq~\la\bfe_2, \gamma_1(t^\sharp, \xi)-\gamma_2(t^\sharp, s^\sharp)\ra~=~0.\enda\right.\eeq
 Computing the $2\times 2$ matrix of partial derivatives at
 the point $ (t^\sharp, s^\sharp, \xi)=(t_0,0,0)$, we obtain
 $${\partial (\Phi_1,\Phi_2)\over \partial (t^\sharp, s^\sharp)}~=~
 \begin{pmatrix}   \la \bfe_1, \bfv_1-\bfv_2\ra && -1\\[3mm]   \la \bfe_2, \bfv_1-\bfv_2\ra
 && 0 \end{pmatrix}.$$
Since $\bfe_2 = \bfn$, by the assumption (\ref{ndeg}) this matrix has full rank.
Hence the functions $\xi\mapsto t^\sharp(\xi)$, $\xi\mapsto s^\sharp(\xi)$ are well
defined.  In turn, by (\ref{gsh})
this determines the endpoint of the new curve $\gamma$ at time
$t^\sharp$.

Next,  conservation of mass and momentum imply
\bel{rsh}
\eta(t^\sharp(\xi), \xi)~=~\eta_1(t^\sharp(\xi),\xi) + \eta_2\bigl(t^\sharp(\xi),
s^\sharp(\xi)\bigr) {d\over d\xi}
s^\sharp(\xi),\eeq
\bel{msh}(\bfv\eta)(t^\sharp(\xi), \xi)~=~(\bfv_1\eta_1)(t^\sharp(\xi),\xi) + (\bfv_2\eta_2)
\bigl(t^\sharp(\xi),s^\sharp(\xi)\bigr) {d\over d\xi}
s^\sharp(\xi).\eeq
We can thus recover the local velocity as a quotient
\bel{vsh}
\bfv(t^\sharp(\xi), \xi)~=~{(\bfv\eta)(t^\sharp(\xi), \xi)\over \eta(t^\sharp(\xi), \xi)}\,.\eeq
Here and in the sequel, $\bfv_1,\bfv_2, \bfv$ refer to the velocities of particles along the singular curves $\gamma_1,\gamma_2,\gamma$.   Similarly,
$\eta_1,\eta_2,\eta$ refer to
the densities of particles along these curves.

\begin{figure}[ht]
\centerline{\hbox{\includegraphics[width=7cm]
{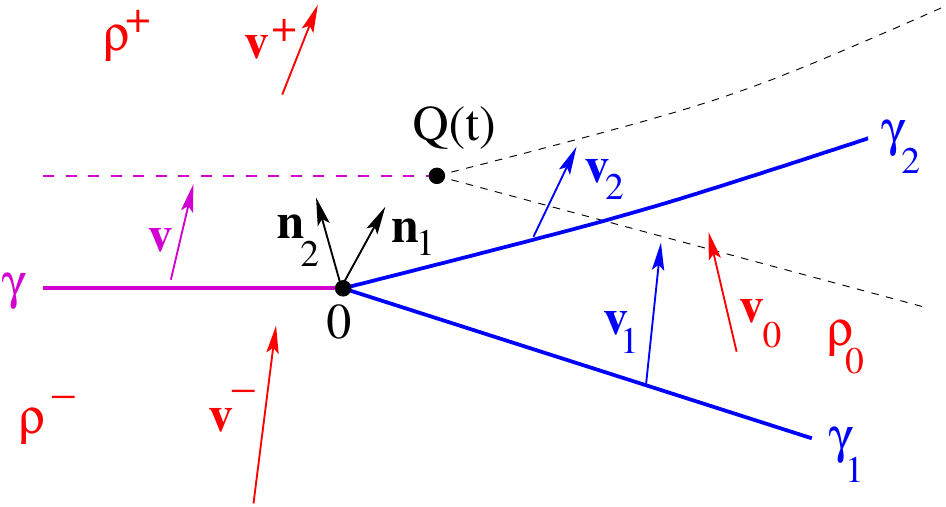}}}
\caption{\small  Two singular curves $\gamma_1,
\gamma_2$, touching at one endpoint, at a nonzero angle. As a result of the collision, more and more particles
from the two curves join together, forming a new singular curve $\gamma$.}
\label{f:gen39}
\end{figure}

\subsection{Two singular curves touching at a positive angle.}

By a change of coordinates we  again assume that, at the initial time $t=t_0$, the curves $\gamma_1,\gamma_2$
can be parameterized as in (\ref{ttc}).    However, instead of (\ref{phi12}), we now assume
\bel{phi3}
\phi_1(0)\,=\,\phi_2(0)\,=\,0, \qquad\qquad \phi_1'(0)\,<\,\phi_2'(0).\eeq
As shown in Fig.~\ref{f:gen39},
we again assume that, in the regions bounded by these two curves, the
density and velocity of the particles are given by analytic functions
$\rho^-, \rho_0, \rho^+$, and $\bfv^-, \bfv_0, \bfv^+$.
This allows us to solve the evolution equations (\ref{ssy}) separately for the two curves.

We denote by $Q(t)$ the intersection point at time $t$.
Calling  $\bfn_1,\bfn_2$
the normal vectors, and $\bfv_1,\bfv_2$ the velocities of particles on the two curves,
the velocity $\dot Q = dQ/dt$
of the intersection point
is determined by the two scalar equations
\bel{dotQ}\bega{rl}\la \dot Q,\,\bfn_1\ra&=~\la \bfv_1,\,\bfn_1\ra,\\[2mm]
\la \dot Q,\,\bfn_2\ra&=~\la \bfv_2,\,\bfn_2\ra.\enda\eeq
The strict inequality in (\ref{phi3}) implies that $\bfn_1\not=\bfn_2$, hence
the vector $\dot Q$ is uniquely determined by (\ref{dotQ}).

We seek conditions on  $\bfn_1,\bfn_2$ and $\bfv_1,\bfv_2$
which guarantee that the maps
\bel{xts}\xi\mapsto t^\sharp(\xi),\qquad\qquad \xi\mapsto s^\sharp(\xi)\eeq
are both strictly increasing.   Notice that this is a key requirement, for the physical admissibility of the solution. Indeed, it implies that as time increases,
particles from both  singular curves join together in the single curve $\gamma$.

Call $\ov \bfn_1,\ov \bfn_2$ the unit normal vectors to the two curves at the origin, at
time $t=t_0$.   Instead of (\ref{Phi12}), we can now determine
the functions $t^\sharp, s^\sharp$ by the equations
\bel{Psi12}\left\{\bega{rl} \Psi_1&\doteq~\la\ov\bfn_1, \gamma_1(t^\sharp, \xi)-\gamma_2(t^\sharp, s^\sharp)\ra~=~0,\\[2mm]
 \Psi_2&\doteq~\la\ov \bfn_2, \gamma_1(t^\sharp, \xi)-\gamma_2(t^\sharp, s^\sharp)\ra~=~0.\enda\right.\eeq
 Computing the $2\times 2$ matrix of partial derivatives at
 the point $ (t^\sharp, s^\sharp, \xi)=(t_0,0,0)$, we obtain
 $${\partial (\Psi_1,\Psi_2)\over \partial (t^\sharp, s^\sharp)}~=~
 \begin{pmatrix}   \la \ov\bfn_1, \bfv_1-\bfv_2\ra && -b\\[3mm]   \la \ov\bfn_2, \bfv_1-\bfv_2\ra
 && 0 \end{pmatrix},$$
where, by (\ref{phi3}),
$$b~=~\left\langle \ov \bfn_1\,,~{\partial \gamma_2(0, s)\over\partial s}
\right\rangle~=~\left\langle {(-\phi_1'(0), 1)\over \sqrt{1+(\phi_1'(0))^2}} \,,~(1, \phi_2'(0))
\right\rangle~=~{\phi_2'(0)-\phi_1'(0)\over \sqrt{1+(\phi_1'(0))^2}}~>~0.$$

 Differentiating (\ref{Psi12}) w.r.t. $\xi$ at the point $ (t^\sharp, s^\sharp, \xi)=(t_0,0,0)$ leads to
\begin{equation*}\left(\begin{array}{cc}\la \ov \bfn_1,\; \bfv_1-\bfv_2\ra&-\left\langle \ov\bfn_1,\;\frac{\partial \gamma_2}{\partial s}\right\rangle\\[2mm]\la \ov\bfn_2,\; \bfv_1-\bfv_2\ra&-\left\langle \ov\bfn_2,\; \frac{\partial \gamma_2}{\partial s} \right\rangle\end{array}\right)\left(\begin{array}{c}\frac{\partial t^{\sharp}}{\partial\xi}\\[2mm]
\frac{\partial s^{\sharp}}{\partial\xi}\end{array}\right)~
=~
-\left(\begin{array}{c}\left\langle \ov\bfn_1,\;\frac{\partial \gamma_1}{\partial\xi} \right\rangle \\[2mm] \left\langle \ov\bfn_2,\; \frac{\partial \gamma_1}{\partial\xi}\right\rangle\end{array}\right).\end{equation*}
Equivalently,
\begin{equation*}\left(\begin{array}{cc}\la \ov \bfn_1,\; \bfv_1-\bfv_2\ra&-b\\[2mm]\la \ov\bfn_2,\; \bfv_1-\bfv_2\ra&0\end{array}\right)\left(\begin{array}{c}\frac{\partial t^{\sharp}}{\partial\xi}\\[2mm]
\frac{\partial s^{\sharp}}{\partial\xi}\end{array}\right)
~=~
-\left(\begin{array}{c}0 \\[2mm] \left\langle \ov\bfn_2,\; \frac{\partial \gamma_1}{\partial\xi}\right\rangle\end{array}\right).\end{equation*}
By  using Cramer's rule, one obtains
\bel{ts}\frac{\partial t^{\sharp}}{\partial\xi}
~=\,-\frac{ \left\langle \ov\bfn_2,\; \frac{\partial \gamma_1}{\partial\xi}\right\rangle}{\la \ov\bfn_2,\; \bfv_1-\bfv_2\ra},\qquad\qquad
\frac{\partial s^{\sharp}}{\partial\xi}
~=\,-\frac{\la \ov\bfn_1,\; \bfv_1-\bfv_2\ra}{ b\,\la \ov\bfn_2,\; \bfv_1-\bfv_2\ra}\left\langle\ov\bfn_2,\; \frac{\partial \gamma_1}{\partial\xi}\right\rangle. \eeq
Recalling (\ref{phi3}),  at the point $(t^{\sharp},s^{\sharp},\xi)=(t_0,0,0)$
we compute
$$\left\langle \ov \bfn_2\,,~{\partial \gamma_1\over\partial \xi}
\right\rangle~=~\left\langle {(-\phi_2'(0), 1)\over \sqrt{1+(\phi_2'(0))^2}} \,,~(1, \phi_1'(0))
\right\rangle~=~-{\phi_2'(0)-\phi_1'(0)\over \sqrt{1+(\phi_2'(0))^2}}~<~0.$$
Since  $b>0$, it follows from (\ref{ts}) that the inequalities   $\frac{\partial t^{\sharp}}{\partial\xi}> 0, \frac{\partial s^{\sharp}}{\partial\xi}> 0$ will hold under the   assumptions
\bel{nv12}\la \ov\bfn_1,\; \bfv_1-\bfv_2\ra~\geq~0,\qquad \la \ov\bfn_2,\; \bfv_1-\bfv_2\ra~>~0.\eeq

\subsection{Local construction of singular curves, near a point of intersection.}
In both of the above cases (i) and (ii), the new singular curve determined by the
interaction can be constructed by solving the system of equations (\ref{ssy}).
However, now the boundary conditions are assigned not at
a fixed initial time $t=t_0$, but along a curve parameterized by
$\xi\mapsto (t^\sharp(\xi), \xi)$.

The local construction of an analytic solution, performed
in Section~\ref{s:3}, can be easily adapted to this somewhat more general situation.

We call
\bel{xn12}x_1\,=\,\xi,\qquad\quad x_2\,=\, \tau\,\doteq\,t-t^\sharp(\xi),\eeq the new
independent variables. Moreover, recalling (\ref{rsh})
and (\ref{vsh}), we
introduce the  new  dependent variables
$$\left\{\bega{rl}Y(t,\xi)&=~y(t,\xi)-\gamma_1(t^\sharp(\xi),\xi),
\\[2mm] W(t,\xi)&=~\eta(t,\xi)-\eta(t^\sharp(\xi),\xi),\\[2mm]
Z(t,\xi)&=~\bfv(t,\xi)-\bfv(t^\sharp(\xi),\xi),\enda\right. \qquad \qquad \left\{\bega{rl}
{\bf P}(t,\xi)&=~\bfp(t,\xi)- y_\xi(t^\sharp(\xi),\xi),\\[3mm] \tau(t,\xi) &= ~t-t^\sharp(\xi)\,.\enda\right.$$
Here it is useful to observe that, knowing the function $\xi\mapsto
y(t^\sharp(\xi),\xi)$, the partial derivative $y_\xi$ is immediately recovered by the formula $$y_\xi(t^\sharp(\xi),\xi)~=~{d\over d\xi}y(t^\sharp(\xi),\xi) - y_t(t^\sharp(\xi),\xi)\frac{d }{d\xi}t^{\sharp}(\xi)
~=~{d\over d\xi}y(t^\sharp(\xi),\xi) - \bfv(t^\sharp(\xi),\xi)\frac{d }{d\xi}t^{\sharp}(\xi).
$$

Starting with  the system (\ref{ssy}), and  adding the trivial equation $\tau_t=1$,
we shall derive a corresponding system of equations for
the vector $\bfU=(Y,W,Z,{\bf P}, \tau)\in \R^8$, in terms of  the independent variables $x_1, x_2 $ in (\ref{xn12}).

Toward this goal, we need to write $\bfp_t=\bfv_{\xi}$ in terms of the variables.  One has
$$({\bf P}+y_{\xi}(t^{\sharp}(\xi),\xi))_t~=~(Z(t,\xi)+\bfv(t^{\sharp}(\xi),\xi))_{\xi},$$
hence
$${\bf P}_t~=~\frac{d}{d\xi}Z(t,\xi)+\frac{d}{d\xi}\bfv(t^{\sharp}(\xi),\xi),$$
where the second term is already known from the boundary conditions.
In turn, this implies
$${\bf P}_{\tau}=Z(\tau+t^{\sharp}(\xi),\xi)_{t}\frac{\partial}{\partial\xi}t^{\sharp}(\xi)+Z(\tau+t^{\sharp}(\xi),\xi)_{\xi}
+\frac{d}{d\xi}\bfv(t^{\sharp}(\xi),\xi)),$$
$${\bf P}_{x_2}=c(x_1) Z_{x_2}+Z_{x_1}+\frac{d}{d\xi}\bfv(t^{\sharp}(\xi),\xi)),\quad c(x_1)=\frac{\partial}{\partial x_1}t^{\sharp}(x_1).$$
As a result, in these new variables  the system (\ref{ssy}) can be written as
$${\bf A}(x_1)\bfU_{x_2}={\bf B} \bfU_{x_1}+{\bf \widetilde{G}}(x_1,\bfU),$$
where ${\bf A}(x_1)$ is a matrix with $a_{ii}=1$ for all $i=1,\cdots,8$, $a_{64}=a_{75}=-c(x_1)$ and all other entries are zero.
Since ${\bf A}(x_1)$ is invertible, we eventually obtain
the following Cauchy problem
\bel{cp4}\left\{\begin{array}{rll}\bfU_{x_2}&=~
{\bf B}(x_1) \bfU_{x_1}+{\bf G}(x_1,\bfU)
\qquad &\hbox{for}~~x_2>0,\\[2mm]\bfU&=~{\bf 0}\qquad
&\hbox{for}~~x_2=0,\end{array}\right.\eeq
where ${\bf B}(x_1)={\bf A}^{-1}(x_1){\bf B}, \;{\bf G}={\bf A}^{-1}(x_1){\bf \widetilde{G}}$.

 A local solution, in the class of analytic functions, is again provided by the
 Cauchy-Kovalevskaya theorem.

\section{Generic singularity formation}
\label{s:5}
\setcounter{equation}{0}

For readers' convenience, we review some basic
definitions and results in transversality theory.
For the proofs, see \cite{Bloom, GG}.

{\bf Definition}. {\it Let $f:\;X\mapsto Y $ be a  smooth map of
manifolds and let $W$ be a submanifold of $Y$. We say that $f$ is
transversal to $W$ at a point $p\in X$, and write $f\pitchfork_p W$,
if
\begin{itemize}\item either $f(p)\notin W$,
\item {\text or else} $f(p)\in W$ and $T_{f(p)}Y=(df)_p(T_pX)+T_{f(p)}W.$\end{itemize}We say that $f$ is transverse to $W$, and
write $f\pitchfork W$, if $f\pitchfork_p W$ for every point $p\in
X$.}

 Here $T_pX$ denotes the tangent space at the
point $p\in X$, while $T_qY$ and $T_qW$ denotes respectively the
tangent spaces to $Y$ and to $W$ at the point $q\in W\subset Y$.
Finally, $(df)_p:\;T_pX\mapsto T_{f(p)}Y$ denotes the differential
of the map $f$ at $p$.

{\bf Transversality Theorem}. {\it  Let $X, \Theta$ and $Y$ be smooth
manifolds, $W$ a submanifold of $Y$. Let $\theta\mapsto
\phi^{\theta}$ be a smooth map which to each $\theta\in \Theta$
associates a function $\phi^{\theta}\in \C^{\infty}(X,Y)$, and
define $\Phi: X\times\Theta\mapsto Y$ by setting
$\Phi(x,\theta)=\phi^{\theta}(x)$. \\
If $\Phi\pitchfork W$, then the
set $\{\theta\in\Theta:\;\phi^{\theta}\pitchfork W\}$ is dense in
$\Theta$.}


 \v
We wish to understand the formation of new singularities for a solution to
(\ref{SP}), for a generic set of
smooth initial data.
Given $T>0$ and a vector field $\bfw\in \C^2_0(\R^2;\R^2)$, for $r=0,1$
we consider the
sets
\bel{Sdef}
\S_r~\doteq~\Big\{(t,x)\,;~~0<t<T,~~\rank (I+ tD\bfw(x))=r\Big\}~\subset~
 \R\times\R^2.\eeq
To analyze these sets, we recall
\begin{lemma}\label{l:corank}
Within the space  of all $2\times 2$ matrices, the set
 \bel{MR}\M_1~\doteq~\{A\in \R^{2\times 2}\,;~~\hbox{\rm rank}(A)=1\}\eeq
 is an embedded submanifold of  dimension 3.
\end{lemma}

We shall apply the above tools from differential geometry
to study the solution to (\ref{SP}), for a generic vector field $\bfw$ in $\R^2$.
Recalling the definitions (\ref{C30}) and (\ref{Sdef}), we begin by proving

\begin{lemma}\label{l:82}  There exists an open dense set of
vector fields $\F_1\subset \C^3_0(\R^2;\R^2)$ such that, for every $\bfw\in\F_1$, one has
\begin{equation}\label{emp}\S_0~\doteq~\Big\{(t,x)\in [0,T]\times \R^2\,;~~\rank\bigl(I+tD_x{\bf w}(x)\bigr)=0\Big\}~=~\emptyset\end{equation}
and, moreover, the set
\bel{SS1}\S_1~\doteq~\Big\{(t,x)\in\, [0,T[\times \R^2\,;~~\rank\bigl(I+tD_x{\bf w}(x)\bigr)=1\Big\}\eeq
is a $\C^2$ embedded manifold of dimension 2.
\end{lemma}

 \v
 {\bf Proof.} {\bf 1.}
 Let a vector field $\bfw\in \C^3_0(\R^2;\R^2)$ be given, and  consider the map
 \bel{em1}(t,x)~\mapsto ~\Phi(t,x)~\doteq~I+t D\bfw(x),\eeq
 from $[0,T]\,\times\R^2$ into $\R^{2\times 2}$.
 We shall construct  a family of perturbed vector fields $\bfv(\theta)$
depending on an auxiliary parameter $\theta\in \R^N$, such that $\bfv(0) = \bfw$
and moreover the map
\bel{trm}(t,x,\theta)~\mapsto ~\Phi(t,x,\theta)~=~I + t [D_x\bfv(\theta)] (x)\eeq
 is transversal to the manifold $\M_0= \{{\bf 0}\}\subset \R^{2\times 2}$
 containing just the zero matrix.

 \v
 {\bf 2.} Consider any point $(t_0,x_0)\in]0,T[\times\R^2$ satisfying
 \bel{wnul}I+t_0D_x{\bf w}(x_0)={\bf 0}.\eeq
 For $\theta=(\theta_{\alpha\beta})\in\R^{2\times2}$, we define a family of perturbed vector fields \bel{per}{\bf v}(\theta)(x)={\bf w}(x)+{\bf p}(x,\theta),\eeq where
 ${\bf p}(\cdot, \theta)\in \C^\infty_c(\R^2;\R^2)$ for every $\theta$,
 and moreover
  \bel{vef}{\bf p}(x ,\theta)~=~\begin{pmatrix} \theta_{11} & \theta_{12}\cr
  \theta_{21} & \theta_{22}\end{pmatrix}\begin{pmatrix} x_1\cr x_2\end{pmatrix}\eeq
 for $x$ in a neighborhood of $x_0$.
   \v
 {\bf 3.} 
 By (\ref{wnul}),
 at the point $(t_0,x_0, \theta)$ one has
 $$I+t[D_x{\bf v}(\theta)](x)~=~I+t_0D_x{\bf w}(x_0)+t_0D_x{\bf p}(x_0,\theta)~=~t_0\begin{pmatrix} \theta_{11} & \theta_{12}\cr
  \theta_{21} & \theta_{22}\end{pmatrix}$$
 Taking the partial derivatives of the map $\Phi=[\Phi_{ij}]$ in
 (\ref{em1}) w.r.t.~the four parameters $\theta_{\alpha \beta}$,
 we obtain
 $$\frac{\partial \Phi_{ij}}{\partial\theta_{\alpha\beta}}~=~\left\{ \bega{cl} t_0\quad
&\hbox{if}\quad i=\alpha
~~\hbox{and}~~j=\beta,\\[2mm]
0\quad &\hbox{otherwise.}\enda\right.$$
 Since $t_0>0$,  we conclude that the map  $\Phi$ at (\ref{trm}) is transversal to the
 zero-dimensional manifold $\M_0$, for $(t,x)$ in a neighborhood $\N_{(t_0,x_0)}$
  of $(t_0,x_0)$ and $\theta$ in a neighborhood of ${\bf 0}\in\mathbb{R}^{2\times2}$.
  \v
  {\bf 4.}  By the definition (\ref{C30}), the set
  $\S_0$ in (\ref{emp})  is compact.
Hence there exist finitely many points $(t_i,x_i)\in \S_0$, $ i=1,2,\cdots,m$, such that
the corresponding open neighborhoods $\mathcal{N}_{(t_i,x_i)}$ cover
$\S_0$.
For each $i$, we construct a $4$-parameter family of
perturbations $\bfp^{(i)}(x,\theta)$ of the vector field $\bfw$, defined as in (\ref{vef}), with
$x_0$ replaced by $x_i$.
Putting together all these perturbations, we thus obtain a family of perturbed vector fields
\bel{vpi}\bfv(\theta)(x)~=~\bfw(x) + \sum_{i=1}^m \bfp^{(i)}(x,\theta)\eeq
depending on
$4m$ parameters.  By construction, for $\theta$ in a neighborhood on ${\bf 0}\in \R^{4m}$,
 the map (\ref{trm})
 is transversal to $\M_0=\{{\bf 0}\}$ as long as  $(t,x)\in\cup_i \N_{(t_i, x_i)}$.
Outside the union of these neighborhoods,
by continuity we have $\Phi(t,x,\theta)\not= {\bf 0}$,
hence transversality trivially holds.

By the Transversality Theorem, for a dense set of values
$\theta\in\R^{4m}$ the map
$(t,x)~\mapsto~\Phi(t,x,\theta)$ is
transversal to the zero-dimensional manifold $\{{\bf 0}\}\subset\R^{2\times 2}$.
Since $[0,T]\times \R^2$  has dimension $3$ and $\R^{2\times2}$ has dimension $4$,
this transversality condition implies
$$\Phi(t,x,\theta )~=~I+t[D_x\bfv(\theta)](x)~\neq~{\bf0}$$
for all $(t,x)\in [0,T]\,\times \R^2$.

We conclude that the set $\F_0$ of all vector fields $\bfw\in\C^3_0$ such that
\bel{q6}
I+tD_x\bfw(x)~\neq~{\bf0}\qquad\forall (t,x)\in [0,T]\times \R^2\eeq
is dense.
\v
{\bf 5.}
To prove that $\F_0$ is open in the topology of $\C^3_0(\R^2;\R^2)$, consider
a  sequence of vector fields
$\{{\bf w}_\nu\}_{\nu\geq1}$ converging to $\bf w$, and assume that
${\bf w}_{\nu}\notin\F_0$ for all $\nu\geq1$. To fix the ideas,
let $(t_{\nu},x_{\nu})\in [0,T]\times \R^2$ be points at which
\begin{equation}\label{condi}I+t_{\nu}D_x{\bf
w}_{\nu}(x_{\nu})={\bf0}.\end{equation}
Since $\bfw_\nu\to \bfw$ in $\C^3_0$ 
it follows that  the convergence $D\bfw_\nu(x)\to 0$ as $|x|\to \infty$ holds
uniformly w.r.t.~$\nu$.  Hence the sequence $(x_\nu)_{\nu\geq 1} $
must be bounded.
By compactness,  by possibly taking a subsequence we can assume
$(t_{\nu},x_{\nu})\rightarrow(\bar{t},\bar{x})$
as $\nu\rightarrow\infty$. By continuity, it follows from (\ref{condi})
that $$I+\bar{t}D_x{\bf
w}(\bar{x})={\bf0}.$$
Hence $\bfw\notin \F_0$.
\v
{\bf 6.} By the previous steps, there exists an open dense set
$\F_0\subset \C^3_0(\R^2;\R^2)$ such that
\bel{jn0}\bfw\in \F_0\qquad\implies\qquad I + tD\bfw(x)\not= {\bf 0}
\qquad\forall (t,x)\in [0,T]\times \R^2.\eeq
In the remainder of the proof we construct an open dense
subset $\F_1\subset \F_0$
such that, for $\bfw\in \F_1$, the set $\S_1$ is a $\C^2$ embedded manifold.

Toward this goal, given $\bfw\in \F_0$,
we will construct  a family of perturbed vector fields $\bfv(\theta)$,
depending on an auxiliary parameter $\theta\in \R^N$,
 such that $\bfv(0) = \bfw$
and moreover
the map
\bel{trv}(t,x,\theta)~\mapsto ~\Phi(t,x,\theta)~=~I + t [D_x\bfv(\theta)] (x)\eeq
 is transversal to the manifold $\M_1$ of all $2\times 2$ matrices with rank 1.
\v
{\bf 7.}  Let $\bfw\in \F_0$ and consider any point $(t_0, x_0) \in \, [0,T]\,\times\R^2$
such that
 \bel{PM0}
\rank\bigl( I+t_0D\bfw(x_0)\bigr)~=~1\,.\eeq
Call
\bel{Ade}
A~=~[a_{ij}]~\doteq~I+t_0D\bfw(x_0).\eeq
Since $A$ is not the zero matrix, one can easily
find another $2\times 2$ matrix $B=[b_{ij}]$  such that
\bel{Btra}\kappa~\doteq~{d\over d\theta} \det(A+\theta B)\bigg|_{\theta=0}~\not= ~0.\eeq
For example, if $a_{12}\not= 0$, it suffices to take $b_{11}=b_{12}=b_{22}=0$,
$b_{21}=1$.

We now
construct a
1-parameter family
of vector fields $\bfp(\cdot,\theta)\in \C^\infty_c(\R^2;\R^2)$ such that
$$\bfp(x,\theta) ~=~\theta B x,$$
for $x$ in a neighborhood of $x_0$.
For $\theta\in \R$ in a neighborhood of the origin,
we then define the perturbed vector fields
$$\bfv(\theta)(x)~=~\bfw(x)+
\bfp(x, \theta).$$
By (\ref{Btra}),  the map
$(t,x,\theta)~\mapsto~I + t[D_x \bfv(\theta)](x)$ is transversal to
the  manifold $\M_1$ at the point $(t_0, x_0,0)$.

By continuity,  the map (\ref{trv}) is still transversal to the
manifold $\M_1$, for $(t,x)$ in a neighborhood $\N_{(t_0, x_0)}$
of $(t_0, x_0)$ and $\theta$ in a neighborhood of  $0$.
\v
{\bf 8.}
Since we are assuming that
$\bfw\in \F_0\subset \C^3_0(\mathbb{R}^2)$, the set
\bel{14}\bega{rl}\ov \S_1&=~\Big\{(t,x)\in [0,T]\times \R^2\,;~~\text{rank}(I+tD_x{\bf w}(x))=1\Big\}\\[3mm]
&=~
\Big\{(t,x)\in [0,T]\times\R^2\,;~~\text{rank}(I+tD_x{\bf w}(x))\leq 1\Big\}
\enda\end{equation}
is compact.  Notice that here we are considering the closed interval $[0,T]$,
instead of the half-open interval used at (\ref{SS1}).
For each $(\bar t,\bar x)\in \ov \S_1$ we can repeat the construction in the previous steps
and find a family of perturbed vector fields $\bfv(\theta)$ depending on
a parameter $\theta\in \R$ such that the
corresponding map in (\ref{trv}) is transversal to $\M_1$ on a neighborhood
of $(\bar t,\bar x)$.

We now choose finitely many points $(t_i, x_i)$,
$i=1,\ldots,\nu$, such that the corresponding
 neighborhoods $\N_{(t_i, x_i)}$ cover the compact set $\ov \S_1$.
 Combining all these perturbations, we thus obtain a family of vector fields
$\bfv(\theta)$ depending
on $\nu$ scalar parameters.   The corresponding map (\ref{trv}) is transversal
to $\M_1$ for all $(t,x, \theta)$ in a neighborhood of $\ov \S_1\times \{{\bf 0}\}$.
On the other hand, for $(t,x)\notin \cup_i \N_i$ and $\theta$ sufficiently
close to ${\bf 0}$, transversality trivially holds.
\v
{\bf 9.} By the Transversality Theorem we conclude that
there exists a dense set ${\bf
\F_1}\subset \F_0\subset \C_0^3(\R^2;\R^2)$ such that, for every ${\bf w}\in\F_1$,
the map (\ref{em1}) is transversal to $\M_1$.

It remains to show that $\F_1$ is open in the topology of $\C^3_0$.
Toward this goal, we first observe that, for any  $\bfw\in \F_0$
(\ref{14}) yields the equivalence
\bel{eq2}
\rank(I + t D_x \bfw(x))~=~1 \qquad\iff\qquad \det  (I + t D_x \bfw(x))~=~0.\eeq
The assumption that $\bfw$ is transversal to $\M_1$
on the domain $[0,T]\times\R^2$  thus becomes equivalent to the
statement that, for every $(t,x)$, the system of  equations
\bel{4e}
\det  (I + t D_x \bfw(x))~=~0, \qquad \partial_t\bigl[\det  (I + t D_x \bfw(x))\bigr]~=~0,
\qquad  \nabla_x\bigl[\det  (I + t D_x \bfw(x))\bigr]~=~0,\eeq
has no solution.
Notice that this condition is precisely what is needed to ensure that the
set $\bigl\{(t,x)\,;~\det  (I + t D_x \bfw(x))=0\bigr\}$ is a 2-dimensional embedded manifold
in $\,]0,T[\times\R^2$.

Consider a convergent sequence $\bfw_\nu\to \bfw$ in $\C^3_0$,
with ${\bf w}_{\nu}\notin\F_1$ for all $\nu\geq1$.
This implies that there exists a sequence of points $(t_\nu, x_\nu)\in [0,T]\times\R^2$
such that, for every $\nu\geq 1$,
\bel{4nu}
\det  (I + t_\nu  D_x \bfw_\nu(x_\nu))~=~0, \quad \partial_t \bigl[\det(I + t_\nu  D_x \bfw_\nu(x_\nu))\bigr]~=~0,
\quad  \nabla_x \bigl[\det (I + t_\nu  D_x \bfw_\nu(x_\nu))\bigr]~=~0.\eeq
Since $\bfw\in \C^3_0$, the first  identities in (\ref{4nu})
imply that the sequence $x_\nu$ is bounded.
By
possibly taking a subsequence we can assume the convergence
$(t_{\nu},x_{\nu})\rightarrow(\bar{t},\bar{x})$, as
$\nu\rightarrow\infty$. By the $\C^3$ convergence $\bfw_\nu\to \bfw$, this implies
that at the point $(t,x)=(\bar t,\bar x)$ all identities in (\ref{4e}) are simultaneously
satisfied.   Hence $\bfw\notin \F_1$ as well.

The above argument shows that $\F_1$ is open, completing the proof.
\endproof
\v

Next, we study the  behavior of a generic solution near a point
$(t_0, x_0 + t_0 \bfw(x_0))$ where a new singularity is formed.
Notice that in this case $t_0$ corresponds to a local minimum of the
map $x\mapsto \tau(x)$, implicitly defined by
the equation $\det(I + \tau D_x\bfw(x))=0$.

\begin{theorem}\label{t:53}
For any given $T>0$ there exists an open dense set $\F\subset\C^3_0(\R^2;\R^2)$
such that, for every vector field $\bfw\in \F$,
the following holds.
\begi
\item[(i)]  The set $\S_1=\bigl\{(t,x)\in [0,T[\,\times\R^2\,;~\rank(I+tD_x{\bf w}(x))=1
\bigr\}$  is an embedded manifold of dimension 2.
\item[(ii)]  There are at most finitely many points $(\tau_k, y_k)\in\S_1$,
$k=1,\ldots,N$ where this manifold is
perpendicular to the $t$-axis.
\item[(iii)]
Call $\tau=\tau(x)$ the function
implicitly defined by $(\tau, x)\in \S_1$ in a neighborhood of
such  points. Then the Hessian matrix of second order partial
derivatives
\bel{HM2}
H~=~\left[\begin{matrix} \tau_{x_1x_1} & \tau_{x_1 x_2}\cr
 \tau_{x_2x_1} & \tau_{x_2 x_2}\end{matrix}\right]\eeq
 has  rank $2$ at each point $y_k$, $k=1,\ldots,N$.
 \item[(iv)] At each point $(\tau_k, y_k)$, $k=1,\ldots,N$, the matrix
$I+\tau_k D_x{\bf w}(y_k)$ has a nonzero eigenvalue.
\endi
\end{theorem}

\v
{\bf Proof.} {\bf 1.} According to  Lemma~\ref{l:82}, there exists an open dense
subset $\F_1\subset\C_0(\R^2;\R^2)$ such that, for every $\bfw\in \F_1$,
the system of four scalar
equations (\ref{4e}) has no solution on $[0,T]\times\R^2$.
By the implicit function theorem,
this implies that $\S_1$ is an embedded manifold.
\v
{\bf 2.} Consider the function
\bel{phidef}\phi(t,x) ~\doteq~\det  (I + t D_x \bfw(x)).\eeq
To prove (ii), we need to show that, for all $\bfw$
in an open dense subset $\F\subset\F_1$, the system of three scalar equations
\bel{5e}
\phi(t,x)~=~\phi_{x_1}(t,x)~=~\phi_{x_2}(t,x)~=~0\eeq
has finitely many solutions on $[0,T]\times\R^2$.

 Given $\bfw\in \F_1$,
consider a point $(t_0, x_0)\in \ov\S_1$ and
the  matrix $A=I+t_0 D_x \bfw(x_0)$. Since $A$ has rank one, we can
write it in the form
\bel{AR1}A~=~\bfb \,\bfe^T,\eeq
where $\bfe,\bfb\in \R^2$ are column vectors, $|\bfe|=1$, and
$^T$ denotes transposition.
For any $x\in \R^2$, this means
$$Ax ~=~\bfb \bfe^Tx~=~\langle \bfe, x\rangle \bfb.$$
Throughout the paper,
$\langle\cdot,\cdot\rangle$ denotes the Euclidean inner product.
Moreover, we shall denote by $\bfb^\perp$ the perpendicular vector, obtained by rotating
$\bfb$ by an angle $\pi/2$, counterclockwise.

We now construct a 3-parameter family of vector fields
$\bfp(\cdot, \theta)\in \C^\infty_c(\R^2;\R^2)$, such that
\bel{166}
\bfp(x,\theta_1,\theta_2,\theta_3)~=~
\theta_1  \langle \bfe^\perp, x-x_0\rangle \bfb^\perp +
\theta_2  \langle\bfe, x-x_0\rangle \langle\bfe^\perp, x-x_0\rangle  \bfb^\perp  +
 \theta_3  \langle\bfe^\perp, x-x_0\rangle ^2 \bfb^\perp,\eeq
for all $x$ in a neighborhood of $x_0$.
In turn, this thus yields a family of perturbations of the vector field $\bfw$
depending on $\theta=(\theta_1,\theta_2,\theta_3)$,
namely
\bel{wp3}\bfv(\theta)~=~
\bfw(x)+\bfp(x,\theta) .\eeq
Setting
$$\phi(t,x,\theta)~\doteq~\det  \bigl(I + t D_x \bfw(x)+ t D_x \bfp(x,\theta)\bigr),$$
we claim that the map
\begin{equation}\label{9}(t,x,\theta)~\mapsto ~\Psi(t,x,\theta)~\doteq~
\bigl(\phi,\,
\phi_{x_1},\,\phi_{x_2}\bigr) \end{equation}
is transversal to the zero-dimensional manifold $\{{\bf 0 }\}\subset
\R^{3}$ on a neighborhood of $(t_0, x_0,0)$.
To prove the claim, we
express the determinant in terms of a cross product:
$$\det(A)~=~A\bfe\times A\bfe^\perp,$$
and use coordinates $(\xi_1,\xi_2)$ corresponding to the orthonormal basis
$\{\bfe,\bfe^\perp\}$.
Computing partial derivatives at
the point $(t,x,\theta)=(t_0, x_0,0)\in [0,T]\times \R^2\times\R^3$,
we now obtain
\begin{equation}\label{20}D_{\theta}\Psi=\left(\begin{array}{ccc} \phi_{\theta_1}& \phi_{\theta_2}&
 \phi_{\theta_3}\vspace{2mm}\\
 \phi_{\xi_1\theta_1}& \phi_{\xi_1\theta_2}& \phi_{\xi_1\theta_3}\vspace{2mm}\\
 \phi_{\xi_2\theta_1}& \phi_{\xi_2\theta_2}& \phi_{\xi_2\theta_3}\vspace{2mm}\\\end{array}\right)~
=~t_0\,|\bfb|^2\cdot\left(\begin{array}{ccc}1&0&0\vspace{2mm}\\
\ast&1&0\vspace{2mm}\\
\ast&\ast&2\vspace{2mm}\\\end{array}\right).\end{equation}
Since this matrix has full rank, we conclude that the map $\Psi$ is
transversal to the zero-dimensional manifold $\{\bf 0\}\subset\R^3$, for $(t,x)$ in a neighborhood
$\N_{(t_0,x_0)}$ of $(t_0, x_0)$ and $\theta$ in a neighborhood of $0$.

\v
{\bf 3.} Consider the  set of all $(t,x)\in [0,T]\times\R^2$
for which all three identities in (\ref{5e}) hold.   By compactness, it can be covered
with finitely many neighborhoods $\N_{(t_i, x_i)}$, $i=1,\ldots,m$.
For each $i$, we construct a family of perturbations $\bfv^i(\theta)$, as in
(\ref{wp3}), such  that the map (\ref{9}) is transversal to the zero manifold
$\{\bf 0\}\subset\R^3$, restricted to $\N_{(t_i, x_i)}$.    As in (\ref{vpi}),
putting together all these
perturbations we obtain a family of vector fields $\bfv(\theta)$ depending on
$m$ parameters, such that the corresponding map (\ref{9}) is transversal
to the zero manifold on the whole domain $[0,T]\times\R^2$.

By the Transversality Theorem, we conclude that
 for a dense set of values
$\theta$ in a neighborhood of the origin in $\mathbb{R}^m$, the map (\ref{9}) is
transversal to the zero-dimensional manifold $\{{\bf
0}\}\subset\mathbb{R}^{3}$.
\v
{\bf 4.} We  now define $\F\subset \F_1$ to be the set of all vector fields
$\bfw\in \F_1$ with the following property:
\begi
\item[{\bf (P)}]  {\it At every $(t,x)\in [0,T]\times\R^2$ where (\ref{5e}) holds,
the Jacobian matrix
$D_{(t,x)}\Psi(t,x)$ has the full rank 3.}
\endi
{}From the previous analysis it follows that $\F$ is dense in $\C^3_0(\R^2;\R^2)$.
Moreover, for every $\bfw\in \F$, the system (\ref{5e}) has at most finitely many
solutions on $[0,T]\times\R^2$.

To prove that $\F$ is open, assume by contradiction that there exists $\bfw\in \F$
and  a sequence of vector fields $\bfw_\nu\notin \F$
converging to $\bfw$ in the $\C^3$ norm.
Then there exist points $(t_\nu, x_\nu)\in [0,T]\times\R^2$
such that the corresponding functions $\Psi_\nu$ in (\ref{5e})
satisfy
$$\Psi_\nu(t_\nu, x_\nu)~=~(0,0,0),\qquad \rank  \bigl[D_{(t,x)}\Psi_\nu (t_\nu, x_\nu) \bigr]\leq 2$$
for all $\nu\geq 1$.   By choosing a subsequence, we can assume the convergence $(t_\nu, x_\nu)\to (\bar t,\bar x)$.
The convergence $\bfw_\nu\to \bfw$ in $\C^3$ now yields
$$\Psi(\bar t,\bar x)~=~(0,0,0),\qquad \rank \bigl[D_{(t,x)}\Psi(\bar t,\bar x)
\bigr]~\leq ~2,$$
against the assumption $\bfw\in \F$.

\v {\bf 5.} To prove  that functions $\bfw\in\F$ also satisfy (iii),
let $x\mapsto \tau(x)$ be the function implicitly defined by
\bel{taud}
\det[ I + \tau D_x\bfw(x)]~=~0.\eeq
Consider  a point $(t,x)$
where
\bel{psi0}\phi(t,x)~=~\phi_{x_1}(t,x)~=~\phi_{x_2}(t,x)~=~0.\eeq
 By construction,  every
function $\bfw\in \F$ satisfies the property {\bf (P)}.
Hence, at every point $(t,x)\in [0,T]\times \R^2$ where (\ref{psi0})
holds,
the matrix
\bel{DPtx}
D_{(t,x)}\Psi~=~\left[\begin{matrix}\phi_t&\phi_{x_1}&\phi_{x_2}\cr
\phi_{x_1 t} & \phi_{x_1 x_1} & \phi_{x_1 x_2}\cr
\phi_{x_2 t}& \phi_{x_2 x_1}& \phi_{x_2 x_2}
\end{matrix}\right]~=~\left[\begin{matrix}\phi_t&0&0\cr
\phi_{x_1 t} & \phi_{x_1 x_1} & \phi_{x_1 x_2}\cr
\phi_{x_2 t}& \phi_{x_2 x_1}& \phi_{x_2 x_2}
\end{matrix}\right] \eeq
has rank 3.
Therefore, the Hessian matrix
$$
\left[\begin{matrix}
 \tau_{x_1 x_1} & \tau_{x_1 x_2}\cr
 \tau_{x_2 x_1}& \ \tau_{x_2 x_2}
\end{matrix}\right]~=~-{1\over\phi_t}\cdot \left[\begin{matrix}
 \phi_{x_1 x_1} & \phi_{x_1 x_2}\cr
 \phi_{x_2 x_1}& \phi_{x_2 x_2}
\end{matrix}\right]$$
has rank 2.
\v
{\bf 6.} Finally, we claim that, by possibly replacing $\F$ with an open dense subset
$\F'\subset\F$, the condition (iv) is also satisfied.
To fix ideas, let $\bfw\in\F$, let $(t_0, x_0)$ be a point where  (\ref{psi0}) holds,
and  call
$$A~=~I+t_0 D_x\bfw(x_0).$$
Since the $2\times 2$ matrix $A$ has rank 1, we can write $A=\bfb \bfe^T$, as in
(\ref{AR1}). We observe that the following conditions are equivalent:
\begi
\item[(i)] $A$ has a nonzero eigenvalue.
\item[(ii)] $A^2$ is not the zero matrix.
\item[(iii)] $\langle \bfe, \bfb\rangle\not= 0$.
\endi
If the vector field $\bfw\in\F$ is such that, at $(t_0, x_0)$ the above conditions fail,
we show it can be approximated by vector fields $\bfv(\theta)$ for which
(i)--(iii) hold.
 Indeed, consider a family of vector fields $\bfp(\cdot, \theta)\in\C^\infty_c(\R^2;\R^2)$ such that
 \bel{pnew}\bfp(x,\theta)~=~ \theta \langle \bfe^\perp, x-x_0\rangle  \bfb+
c_1 (\theta) \langle\bfe, x-x_0\rangle \langle\bfe^\perp, x-x_0\rangle  \bfb^\perp  +
c_2 (\theta) \langle\bfe^\perp, x-x_0\rangle ^2 \bfb^\perp,\eeq
for $x$ in a neighborhood of $x_0$.  Define
$$\bfv(\theta)~ \doteq ~\bfw+  \bfp(x,\theta).$$
For a suitable choice of the functions
$c_1, c_2$,
we claim that, replacing $\bfw$ with $\bfv(\theta)$ for
$\theta\approx 0$, at the point $(t_0, x_0)$ the conditions
$\phi=\phi_{x_1}=\phi_{x_2}=0$ still hold.   Moreover, when $\theta\not= 0$  the matrix
$$A(\theta)~=~I+t_0 D_x \bfw(x_0) +  t_0 D_x \bfp(x_0,\theta) $$
satisfies  the equivalent properties  (i)--(iii).

Indeed, for every $\theta$ we still have
$$\hbox{range}(A(\theta)) ~=~\hbox{span} \{\bfb\}$$
However, for $\theta\not= 0$ at the point $x_0$ we have
$$A(\theta) ~=~\bfb \,\bfe^T + \bfb(t_0\theta\bfe^\perp)^T.$$
By assumption, $\langle \bfe,\bfb\rangle=0$. Therefore
$$\langle \bfe + t_0\theta\bfe^\perp,\, \bfb\rangle~=~t_0 \theta
\langle\bfe^\perp,\, \bfb\rangle~=~ \pm t_0\theta|\bfb|
~\not=~0.$$
This implies that $A$ has a nonzero eigenvalue.

Notice that in the above computation the functions $c_1,c_2$ do not play any role.
It now remains to show that, by a suitable choice of $c_1(\theta),c_2(\theta)$,
at the point $(t_0, x_0)$ the identities
 $\phi_{x_1}=\phi_{x_2}=0$
also hold. Toward this goal we compute
$${d\over dh} \phi(t_0, x_0 + h\bfe),\qquad \qquad {d\over dh} \phi(t_0, x_0 + h\bfe^\perp),$$
and show that these derivatives can be made to be zero  for all $\theta\approx 0$ at $h=0$,
by a suitable choice of the functions  $c_1, c_2$.
Differentiating (\ref{pnew}) one obtains
$$\bega{rl}D_x\bfp(x,\theta)& =~\theta\bfb (\bfe^\perp)^T+c_1 (\theta)\langle\bfe^\perp, x-x_0\rangle  \bfb^\perp \bfe^T\\[2mm]
&\qquad +c_1 (\theta) \langle\bfe, x-x_0\rangle  \bfb^\perp(\bfe^\perp)^T+2c_2 (\theta) \langle\bfe^\perp, x-x_0\rangle  \bfb^\perp(\bfe^\perp)^T.\enda$$
Therefore
 $$D_x\bfp(x_0+h\bfe,\theta)~=~\theta\bfb (\bfe^\perp)^T+hc_1(\theta)\bfb^\perp(\bfe^\perp)^T,\quad $$
 $$D_x\bfp(x_0+h\bfe^\perp,\theta)~=~\theta\bfb (\bfe^\perp)^T
 +h\bfb^\perp\Big(c_1(\theta)\bfe^T+2c_2(\theta)(\bfe^\perp)^T\Big).$$
At $h=0$ we thus have
$$\bega{l}\ds {d\over dh} \phi(t_0, x_0 + h\bfe)~=~{d\over dh}A(\theta)\,
\bfe\times A(\theta)\,\bfe^\perp+A(\theta)\,\bfe\times {d\over dh}A(\theta)\,\bfe^\perp
\\[3mm]
 \qquad \ds=~{d\over dh}A_0\,\bfe\times\Big(A_0\,\bfe^\perp+\theta t_0\bfb\Big)
 +A_0\,\bfe\times\Big({d\over dh}A_0\,\bfe^\perp+t_0c_1(\theta)\bfb^\perp\Big)
 \\[3mm]
 \qquad\ds  =~\left({d\over dh}A_0\,\bfe\times A_0\,\bfe^\perp+A_0\,\bfe\times{d\over dh}A_0\,\bfe^\perp\right)+
 t_0\left(\theta{d\over dh}A_0\,\bfe\times\bfb+c_1(\theta)\bfb\times\bfb^\perp\right),
 \enda $$
where $A_0\doteq I+t_0D_x\bfw(x_0+h\bfe)$.
 By assumption, $${d\over dh}A_0\,\bfe\times A_0\,\bfe^\perp+A_0\,\bfe\times{d\over dh}A_0\,\bfe^\perp~=~0.$$
Hence, at $h=0$ we have
 \begin{eqnarray*}{d\over dh} \phi(t_0, x_0 + h\bfe)&=&t_0\left(\theta{d\over dh}A_0\,\bfe\times\bfb+c_1(\theta)\bfb\times\bfb^\perp\right) \\
 &=&t_0\left(\theta t_0\bfe^TD_x^2\bfw(x_0)\,\bfe\times \bfb+c_1(\theta)\bfb\times \bfb^\perp\right)\\
 &=&0,\end{eqnarray*}
provided that we choose
$$c_1(\theta)=-\frac{ t_0\bfe^TD_x^2\bfw(x_0)\,\bfe\times \bfb}{\bfb\times \bfb^\perp}\,\theta.$$
Note that this is meaningful because $\bfb\times \bfb^\perp\neq0$.

A similar calculation yields, always at $h=0$,
\begin{eqnarray*}{d\over dh} \phi(t_0, x_0 + h\bfe^\perp)&=&t_0\left(\theta{d\over dh}A_0\,\bfe\times\bfb+2c_2(\theta)\bfb\times\bfb^\perp\right) \\
 &=&t_0\left(\theta t_0\bfe^TD_x^2\bfw(x_0)\,\bfe\times \bfb+2c_2(\theta)\bfb\times \bfb^\perp\right)\\
 &=&0,\end{eqnarray*}
provided that
  $$c_2(\theta)=-\frac{ t_0\bfe^TD_x^2\bfw(x_0)\,\bfe\times \bfb}{2\,\bfb\times \bfb^\perp}\,\theta.$$

By the previous analysis, given $\bfw\in\F$, for any small $\theta\not= 0$
 the perturbed vector field $\bfv(\theta)=\bfw+\theta \bfp(\cdot, \theta)$
satisfies the additional property (iv).
It is now immediate to check that the set of vector fields $\bfv\in\F$ that
satisfy the additional property (iv) is open w.r.t.~the $\C^3$ norm.
This completes the proof of the Theorem.
\endproof

{\bf Acknowledgments.} The research by the first author
 was partially supported by NSF with
grant  DMS-2006884, ``Singularities and error bounds for hyperbolic equations". The work by the second author was partially supported by NSF with grants DMS-2008504 and DMS-2306258. The work of third author was in part supported by Zhejiang Normal University with grants YS304222929 and ZZ323205020522016004.
\v

\end{document}